\documentclass[12pt]{amsart}
\usepackage[utf8]{inputenc}
\usepackage{enumerate,enumitem,amsmath,amssymb,mathrsfs,stmaryrd}
\usepackage{amsthm}
\usepackage{latexsym,epsf,graphicx,comment,appendix}

\usepackage{hyperref} % Optional: for clickable links
\usepackage{xcolor}

\hypersetup{
        colorlinks=true, % Colors the links instead of putting boxes around them
        linktoc=all,     % Makes both section title and page number clickable in TOC
        linkcolor=blue,  % Color of internal links
        citecolor=blue, % Color of citation links
        filecolor=magenta, % Color of file links
        urlcolor=black    % Color of URL links
    }

\newfont{\bb}{msbm10 at 12pt}
\newfont{\tbb}{msbm10 at 8pt}

\def\R{\r}
\def\H{\mathbb{H}}

\def\R{\mathbb{R}}
\def\G{\mathbb{G}}
\def\L{\mathcal{L}}
\def\C{\mathcal{C}}

\baselineskip 
\usepackage {amsmath}
\usepackage {amsthm}
\usepackage{times}
\usepackage{amscd}
\usepackage{epsf}

\numberwithin{equation} {section}

\begin{document}
\mbox{}\vspace{0.2cm}\mbox{}

\providecommand{\keywords}[1]
{
  \small	
  \textbf{\textit{Keywords---}} #1
}

\theoremstyle{plain}\newtheorem{lemma}{Lemma}[section]
\theoremstyle{plain}\newtheorem{proposition}{Proposition}[section]
\theoremstyle{plain}\newtheorem{theorem}{Theorem}[section]

\theoremstyle{plain}\newtheorem*{theorem*}{Theorem}
\theoremstyle{plain}\newtheorem*{main theorem}{Main Theorem} 
\theoremstyle{plain}\newtheorem*{lemma*}{Lemma}
\theoremstyle{plain}\newtheorem*{remark*}{Remark}
\theoremstyle{plain}\newtheorem*{claim}{Claim}
\theoremstyle{plain}\newtheorem{corollary}{Corollary}[section]
\theoremstyle{plain}\newtheorem*{corollary-A}{Corollary A}
\theoremstyle{plain}\newtheorem*{corollary-B}{Corollary B}

\theoremstyle{plain}\newtheorem*{theorem 2.2}{Theorem 3.1}
\theoremstyle{plain}\newtheorem*{theorem 2.3}{Theorem 3.2}
\theoremstyle{plain}\newtheorem*{theorem 4.2}{Theorem 4.2}
\theoremstyle{plain}\newtheorem*{theorem 5.4}{Theorem 5.1}
\theoremstyle{plain}\newtheorem*{theorem 5.2}{Theorem 5.2}
\theoremstyle{plain}\newtheorem*{theorem 5.3}{Theorem 5.3}

\theoremstyle{plain}\newtheorem*{corollary 2.1}{Corollary 2.1}
\theoremstyle{plain}\newtheorem*{corollary 3.1}{Theorem 4.1}

\theoremstyle{definition}\newtheorem{example}{Example}[section]
\theoremstyle{definition}\newtheorem{remark}{Remark}[section]
\theoremstyle{definition}\newtheorem{definition}{Definition}[section]
\theoremstyle{definition}\newtheorem{acknowledge}{Acknowledgment}
\theoremstyle{definition}\newtheorem{conjecture}{Conjecture}

\begin{center}
\rule{15cm}{1.5pt} \vspace{.3cm}

{\bf\Large Remarks on thin sets in potential theory for sub-Laplacians} 
\vskip .2cm

Shiguang Ma$\mbox{}^\dag$ and Jie Qing$\mbox{}^\ddag$
\vspace{0.2cm} 
\rule{15cm}{1.5pt}
\end{center}

%\acknowledgment{}

\title{}

\begin{abstract} In contrast to Euclidean cases, we observe in this note that,
there are sets in Heisenberg groups that are thin at the origin and block all geodesics from 
the origin. However, in polarizable Carnot groups, using horizontal radial 
curves constructed for polar coordinates by Balogh and Tyson, 
we show that sets, that are thin at the origin, cannot block all horizontal curves. Consequently, 
we are able to use the asymptotic behavior with respect to the fine topology to estimate the Hausdorff
dimension of singularities of potentials for sub-Laplacians and obtain geometric applications. 
\end{abstract}


\maketitle

%\tableofcontents

%%%%%%%%%%%%%%%%%%%%%%%%%%%%%%%%%%%%%%%%%%%%%
%%%%%%%%%%%%%%%%%%%%%%%%%%%%%%%%%%%%%%%%%%%%%

\section{Introduction}\label{Sect:intro}

In this note we want to study thin sets in potential theory for sub-Laplacians in Carnot groups. Our main
reference on potentials for sub-Laplacians in homogeneous Carnot groups is \cite{BLU}. Our study results in applications in geometries modeled 
after those on the boundary at infinity of rank one symmetric spaces.   

Let $\G = (\R^N, \circ, \delta_\lambda)$ be a homogeneous Carnot group, and let $\Gamma$ be the fundamental 
solution to a sub-Laplacian $\L$
on $\G$. For a Radon measure $\mu$ with the support in a domain $D\subset \R^N$, the potential of $\mu$ 
for the sub-Laplacian $\L$ is
$$
R^{\mu, D}_\L (x) = \int_D \Gamma(y^{-1}\circ x) d\mu(y).
$$
We would like to make the note that we always assume through out this paper that Radon measures we consider are such that their potentials are not identically infinite.   

From \cite[Definition 12.1.1]{BLU}, the $\L$-fine topology is the smallest topology on $\R^N$
with respect to which all potentials for sub-Laplacian $\L$ are continuous. A set $E$ is $\L$-thin at a point $x_0$ 
if $x_0$ is not a limit
point of $E$ with respect to the $\L$-fine topology according to \cite[Definition 12.2.1]{BLU}. The Wiener 
criterion for sub-Laplacians \cite[Theorem 12.4.2]{BLU}
is fundamental and characterizes thin sets in terms of $\L$-capacity in Carnot groups. 

We start with the study of the singular behavior of potentials for sub-Laplacians with respect to the $\L$-fine topology. 
The first one is the extension of \cite[Theorem 5.3 in Chapter 2]{Mi}. But, to us, it is originated from \cite{AH} on subharmonic functions
in 2 dimensions.

\begin{theorem 2.2} Let $R^{\mu, D}_\L$ be the potential of a Radon measure $\mu$ with the support in $D$ for the sub-Laplacian $\L$ 
in a homogeneous Carnot group $\G = (\R^N, \circ, \delta_\lambda)$. Then, for any $x_0\in D$, there exists a set $E$ that is $\L$-thin 
at $x_0$ such that 
\begin{equation}\label{Equ:intr-2.2}
\lim_{x\to x_0, x\notin E} R^{\mu, D}_\L(x)/\Gamma (x_0^{-1}\circ x) = \mu(\{x_0\}).
\end{equation}
\end{theorem 2.2}

The second one is a strong improvement of \eqref{Equ:intr-2.2} with the help of a measure theoretic argument (see Lemma \ref{Lem:kpata-1.4}).
It is the extension of \cite[Theorem 1.1]{MQ-2025}, which was responsible for estimating of the size of singularities in higher dimensions.

\begin{theorem 2.3}
Let $R^{\mu, D}_\L$ be the potential of a Radon measure $\mu$ with the support in $D$ for the sub-Laplacian $\L$ in a homogeneous
Carnot group $\G = (\R^N, \circ, \delta_\lambda)$.  Let $S$ be a compact subset in $D$ whose Hausdorff dimension  
is greater than $s$ for $s\in (0, Q-2)$, where $Q$ is the homogeneous dimension of $\G$.  
Then there is a point $x_0\in S$ and a set $E$ that is $\L$-thin at $x_0$ such that 
\begin{equation}\label{Equ:intr-2.3}
R^{\mu, D}_\L (x) \leq  C(d_\L(x_0^{-1}\circ x))^{2 - Q + s}
\end{equation}
for $x$ close to $x_0$ and not in $E$, where $d_\L$ is the $\L$-gauge in \eqref{Equ:L-gauge}.
\end{theorem 2.3}

The Hausdorff measures and the Hausdorff dimensions in Theorem \ref{Thm:theorem-7} are with respect to the $\L$-gauge 
$d_\L$ (cf. \cite[Chapter 13]{BLU}). A result about potentials on Heisenberg groups similar to Theorem \ref{Thm:theorem-7} 
was established before in \cite[Lemma 1.10]{CY}. From the geodesic dynamics of Heisenberg groups $\H^n$ described in \cite{RM}, one observes the following

\begin{corollary 3.1} In Heisenberg group $\mathbb{H}^n$, the set $E$ in \eqref{Equ:set-E}
is $\Delta_\H$-thin at the origin according to Definition \ref{Def:thinness}. However, any geodesic from the origin intersects with $E$.
\end{corollary 3.1}

This is in contrast to the existence of rays escaping from a thin set in Euclidean spaces
\cite[Lemma 8.2 in Chapter 2]{Mi} and \cite[Theorem 1.2]{MQ-2025}, which is
the critical property of thin sets first used in \cite[Theorem 1.1]{AH} in the study of the singular behavior of subharmonic functions
that led to the simplified proof of Huber's theorem on open surfaces \cite{Hu} and later empowered \cite{MQ-2025, MQ-survey} to 
develop applications in conformal geometry in general dimensions. 

On the existence of horizontal curves (not necessarily geodesics) escaping from a thin set, one realizes that polar coordinates of horizontal curves, 
if exist, are invaluable. Indeed the polar coordinates constructed by the horizontal flow lines \eqref{Equ:flow-BT} in \cite{BT}
enable us to establish the existence of horizontal curves escaping from a thin set in polarizable Carnot groups.

\begin{theorem 5.4} Let $\G = (\R^N, \circ, \delta_\lambda)$ be a polarizable Carnot group. And let $E$ 
be a set that is thin at the origin in $\R^N$. Then there always exist horizontal radial curves which are flow lines given by
\eqref{Equ:flow-BT} and do not intersect with $E$, at least in a small neighborhood of the origin.
\end{theorem 5.4}

Consequently we are able to use Theorem \ref{Thm:theorem-7} to estimate the Hausdorff
dimension of singularities of potentials for sub-Laplacians in polarizable Carnot groups (cf. Theorem \ref{Thm:main-potential}) 
using the same ideas in \cite{AH, MQ-2025, MQ-survey, CY}. Notice that the nilpotent part of Iwasawa decomposition of
rank one simple groups is H-type (cf. \cite{CDKR-91}) and therefore polarizable due to \cite[Proposition 5.6]{BT}. 
Thus this gives an alternative proof of the main geometric results in \cite{MQ-sub}. 
For background about scalar curvature and Yamabe type equations in geometries modeled after 
those on the boundary at infinity of rank one symmetric spaces, readers are referred to \cite[Section 4]{MQ-sub} and references therein.

\begin{theorem 5.3} (\cite[Theorem 4.2]{MQ-sub}) Assume $\mathsf{N}$ 
is the nilpotent part of the Iwasawa decomposition of rank one groups:
$SU(n+1, 1)$, $Sp(n+1, 1)$ and $\mathbb{F}_{4(-20)}$. 
Let $S$ be a compact subset of a bounded and open subset $\Omega$ in $\mathsf{N}$. 
And let $\theta = u^\frac 4{Q-2}\theta_0$ be a 
conformal contact form on $\Omega\setminus S$, where $\theta_0$ is the standard contact form on $\mathsf{N}$ and $Q$ is the homogeneous 
dimension of $\mathsf{N}$. 
Suppose that the sub-Riemannian metric associated with $\theta$ is 
geodesically complete near $S$. Then
\begin{equation}\label{Equ:intr-dim-CR}
dim_{\mathcal{H}_{\Delta_{\mathsf{N}}}} (S) \leq \frac {Q-2}2
\end{equation}
when the scalar curvature $R_\theta$ of the conformal contact form $\theta$ is nonnegative.
\end{theorem 5.3}

Theorem \ref{Thm:main-geometry} in the cases of CR geometry has appeared in \cite[Theorem 1.1]{CY}.
The organization of this note is as follows: In Section 2, we recall definitions in potential theory for sub-Laplacians 
from \cite{BLU}. In Section 3, we derive estimates on potentials for sub-Laplacians with respect to the fine topology similar to
those on Riesz potentials in Euclidean spaces. In Section 4, we use the presentation of geodesics from the
origin in Heisenberg groups in \cite{RM} and construct a set that is $\Delta_\H$-thin at the origin and blocks all geodesics from
the origin. In Section 5, we introduce horizontal radial curves for polar coordinates constructed by Balogh and Tyson
in \cite{BT} and use Theorem \ref{Thm:capacity-estimate} to derive consequences in potential theory as well as in geometry.

%%%%%%%%%%%%%%%%%%%%%%%%%%%%%%%%%%%% 
%%%%%%%%%%%%%%%%%%%%%%%%%%%%%%%%%%%%

\section{potentials for sub-Laplacians in homogeneous Carnot groups} \label{Sect:h-carnot}

In this section we mostly follow the book \cite{BLU} to quickly introduce the potential theory for sub-Laplacians in homogeneous Carnot groups. 

%%%%%%%%%%%%%%%%%%%%%%%%%%%%%%%%%%%%%%%%%%%%%%%%%%%%%%

\subsection{Basic notations and definitions} \label{SubSect:definition} Let us start with the definition of homogeneous Carnot groups.

\begin{definition} \label{Def:hcg} (\cite[Definition 1.4.1]{BLU}) We say a Lie group on $\mathbb{R}^N$, $\G= (\mathbb{R}^N, \  \circ, \delta_\lambda)$, 
is a homogeneous Carnot group if the following properties hold
\begin{enumerate}
\item $\mathbb{R}^N$ can be split as $\mathbb{R}^{N_1}\times\mathbb{R}^{N_2}\cdots\times\mathbb{R}^{N_r}$ and the dilation
$$
\delta_\lambda (x^{(1)}, x^{(2)}, \cdots x^{(r)}) = (\lambda x^{(1)}, \lambda^2 x^{(2)}, \cdots \lambda^r x^{(r)}), \text{ where $x^{(i)}\in \mathbb{R}^{N_i}$,}
$$
is an automorphism of the group $\mathbb{G}$ for all $\lambda > 0$. Hence the group with this dilation $\delta_\lambda$ 
is a homogeneous Lie group (cf. \cite[Definition 1.3.1]{BLU}). 
\item Let $Z_1, Z_2, \cdots, Z_{N_1}$ be the left invariant vector fields on $\mathbb{G}$ such that
$$
Z_i (0) = \frac {\partial}{\partial x_i}|_{0} \text{ for $i= 1, 2, \cdots N_1$}.
$$
Then the rank of the Lie algebra generated by $Z_1, Z_2, \cdots Z_{N_1}$ is $N$. 
\end{enumerate}
The group together with the dilation $(\mathbb{R}^N, \circ, \delta_\lambda)$ is said to be a homogeneous Carnot group of step $r$ and $N_1$ generators. And
the vector fields $Z_1, Z_2, \cdots Z_{N_1}$ are called the Jacobian generators.
\end{definition}

By \cite[Proposition 1.3.21]{BLU}, the Lebesgue measure on $\mathbb{R}^N$ is the Haar measure 
on a homogeneous Carnot group $\mathbb{G} = (\mathbb{R}^N, \circ, \delta_\lambda)$. And, if we denote by $|E|$ the Lebesgue measure of a measurable set $E\subset \mathbb{R}^N$,
then
$$
|\delta_\lambda (E)| = \lambda^Q|E|,
$$
where 
$Q= \sum_{i=1}^r iN_i$
is called the homogeneous dimension of $\mathbb{G}$. For our purpose we assume $Q>2$ throughout this note. 
Next important objects for us are sub-Laplacians.  

\begin{definition}\label{Def:sub-L} (\cite[Definition 1.5.1]{BLU})  Let $Z_1, Z_2, \cdots, Z_{N_1}$ be the Jacobian generators of a homogeneous Carnot group
$\mathbb{G} = (\mathbb{R}^N, \circ, \delta_\lambda)$. The operator $\Delta_{\mathbb{G}} = \sum_{i=1}^{N_1} Z_i^2$
is called the canonical sub-Laplacian on $\mathbb{G}$. The operator $\mathcal{L} = \sum_{i=1}^{N_1} Y_i^2,$
where $\{Y_1, Y_2, \cdots, Y_{N_1}\}$ is a frame of left-invariant vector fields in the horizontal distribution and 
generate the Lie algebra of rank $N$, is said to be a sub-Laplacian on $\mathbb{G}$.
\end{definition}

Next are fundamental solutions and gauge functions. 

\begin{definition}\label{Def:fund-sol} (\cite[Definition 5.3.1]{BLU})
Let $\mathbb{G} = (\mathbb{R}^N, \circ, \delta_\lambda)$ be a homogeneous Carnot group. And let $\mathcal{L}$ be a sub-Laplacian 
on $\mathbb{G}$. A function
$\Gamma: \mathbb{R}^N\setminus \{0\}\to \mathbb{R}$ is said to be a fundamental solution for $\mathcal{L}$ if:
\begin{itemize}
\item $\Gamma$ is smooth away from the origin;
\item $\Gamma \in L^1_{\text{loc}}(\mathbb{R}^N)$ and $\Gamma (x)\to 0$ as $x\to \infty$;
\item $- \mathcal{L} \Gamma = \text{Delta}_0,$ where $\text{Delta}_0$ is the Dirac measure supported at the origin $0$. More explicitly,
$\int_{\mathbb{R}^N} \Gamma \mathcal{L}\phi = -\phi(0), \ \forall \phi\in C^\infty_0(\mathbb{R}^N)$.
\end{itemize}
\end{definition} 

It turns out, by \cite[Theorem 5.3.2 and Proposition 5.3.10]{BLU}, for a sub-Laplacian $\mathcal{L}$, 
there always exists a unique fundamental solution, which is positive with the homogeneity 
\begin{equation}\label{Equ:Gamma-scale}
\Gamma (\delta_\lambda (x)) = \lambda^{2 - Q}\Gamma (x)
\end{equation}
(cf. \cite[Proposition 5.3.12 and Proposition 5.3.13]{BLU}). 
In fact, for a sub-Laplacian $\mathcal{L}$,
\begin{equation}\label{Equ:L-gauge}
d_{\L} (x, y)= d_\L (y^{-1}\circ x) = \left\{\aligned \Gamma(y^{-1}\circ x)^{\frac 1{2-Q}} & \text{ when $x\neq 0$}\\
0 \quad & \text{ when $x=y$}\endaligned\right.
\end{equation}
is an $\mathcal{L}$-gauge (cf. \cite[Definition 5.4.1]{BLU}), which is a symmetric homogeneous norm on $\mathbb{G}$. It is also good to know that, 
by \cite[Proposition 5.3.10]{BLU}, $\mathcal{L}$-gauges are unique up to a constant multiple. 

At this point we are ready to introduce potentials of Radon measures for sub-Laplacians in homogeneous Carnot groups.

\begin{definition}\label{Def:potential} Let $\mathbb{G}= (\mathbb{R}^N, \circ, \delta_\lambda)$ be a 
homogeneous Carnot group. And let $\L$ be a sub-Laplacian on $\mathbb{G}$ and $\Gamma$ is its fundamental solution. We define the potential
for the sub-Laplacian $\L$ as
\begin{equation}\label{Equ:potential}
R_\L^{\mu, D} (x) = \int_D \Gamma (y^{-1}\circ x)d\mu(y)
\end{equation}
for all $x\in \R^N$, where $\mu$ is a Radon measure and  $D$ is the support of $\mu$. 
\end{definition}

It is important to know by \cite[Theorem 9.4.4]{BLU} that an $\L$-superharmonic function (cf. (\cite[Definition 7.2.2 and Remark 7.2.3]{BLU}) 
is represented by a $\L$-potential $R^{\mu, D}_\L$ for a nonnegative Radon measure $\mu$ up to $\L$-harmonic functions.  

%%%%%%%%%%%%%%%%%%%%%%%%%%%%%%%%%%%%%%%%%%%%%%%%%%%%%%%

\subsection{$\L$-capacity, $\L$-thin sets, and the $\L$-fine topology}\label{Subsect:cap-thin} 

We refer readers to \cite[Chapter 11 and 12]{BLU} for a more comprehensive introduction and discussion of the $\L$-capacity, $\L$-thinness, and the $\L$-fine topology 
for a sub-Laplacian $\L$ in $\G= (\R^N, \circ, \delta_\lambda)$. The $\L$-fine topology is characterized as the smallest
topology with respect to which all $\L$-superharmonic functions is continuous (cf. \cite[Definition 12.1.1]{BLU}). 
$\L$-thinness is introduced by the $\L$-fine topology (cf. \cite[Definition 12.2.1]{BLU}). It is deep and powerful in the potential theory that 
the Wiener criterion for sub-Laplacians provides an equivalent definition of $\L$-thinness in terms of $\L$-capacity (cf. \cite[Theorem 12.4.2]{BLU}).   
For the introduction of exterior $\L$-capacity for subset in $\R^N$, we refer readers to \cite[Definition 11.5.4 and Definition 11.5.14]{BLU}. For 
simplicity, we summarize it as follows.

\begin{definition}\label{Def:l-capacity} (\cite[Chapter 11]{BLU}) Let $\mathbb{G}= (\mathbb{R}^N, \circ, \delta_\lambda)$ be a 
homogeneous Carnot group. And let $\L$ be a sub-Laplacian on $\mathbb{G}$. First, for any compact subset $K\subset \R^N$, we define the $\L$-capacity as
$$
cap_\L(K) = \inf\{ \mu(K): \mu \text { is a nonnegative Radon measure on $K$ and } R_\L^{\mu, K}\geq 1 \text{ on } K\}
$$
or equivalently (cf. \cite[Theorem 11.5.6]{BLU})
$$
cap_\L(K) = \sup\{ \mu(K): \mu \text { is a nonnegative Radon measure on $K$ and } R_\L^{\mu, K}\leq 1 \text{ in } \R^N\}.
$$
Then, for any subset $E\subset \R^N$, we define the exterior $\L$-capacity as 
$$
cap_\L^* (E) = \inf_{U\supset E \text{ and open }} \, \sup_{K\subset U \text{ and compact}} cap_\L (K).
$$
\end{definition}

To set notations for introducing $\L$-thin sets via the exterior $\L$-capacity, for $x_0\in \R^N$, we recall from \cite[(12.16)]{BLU},
\begin{equation}\label{Equ:gauge-dyadic} 
\C_n  = \{ y\in\R^N: \alpha^n \leq \Gamma (y^{-1}\circ x_0) \leq \alpha^{n+1}\} =  \{ y\in\R^N: \bar \alpha^{-n-1} \leq d_\L(x_0, y) \leq \bar \alpha^{-n}\}
\end{equation} 
for 
\begin{equation}\label{Equ:alpha-bar}
\bar \alpha = \alpha^{\frac 1{Q-2}} > \beta,
\end{equation}
where $\beta (\geq 1)$ is the constant in the improved pseudo-triangle inequality for the $\L$-gauge $d_\L$ of \eqref{Equ:L-gauge} 
(cf. \cite[Proposition 5.14.1 and (12.15a)]{BLU}). 
We remark here that $\bar \alpha$ was taken as $2$ and $\C_n$ are the standard dyadic spherical shells in Euclidean cases. 

\begin{definition}\label{Def:thinness} (\cite[Definition 12.2.1 and Theorem 12.4.2]{BLU}) Let $\mathbb{G}= (\mathbb{R}^N, \circ, \delta_\lambda)$ be a 
homogeneous Carnot group. And let $\L$ be a sub-Laplacian on $\mathbb{G}$. Then we say a subset $E\subset\R^N$ is $\L$-thin at a point $x_0\in \R^N$
if
\begin{equation}\label{Equ:def-l-thin-gauge}
\sum_{n=1}^\infty \alpha^n cap_\L^*(E\cap \C_n) = \sum_{n=1}^\infty \bar\alpha^{n(Q-2)} cap_\L^*(E\cap \C_n)< \infty
\end{equation}
where $\C_n$ given in \eqref{Equ:gauge-dyadic} is the dyadic spherical shell with respect to $d_\L$ and centered at $x_0$.
\end{definition} 

%%%%%%%%%%%%%%%%%%%%%%%%%%%%%%%%%%%%%%%%%%%%%%
%%%%%%%%%%%%%%%%%%%%%%%%%%%%%%%%%%%%%%%%%%%%%%

\section{Asymptotic behavior of potentials with respect to the fine topology}\label{Sect:potential-estimates}

The notion of the fine topology is such an sharp way to study the continuity of potentials due to the well-known
Wiener criterion (cf. \cite[Chapter 12]{BLU} and \cite[Chapter 2]{Mi}).
It turns out the fine topology is also sharp in describing the singular behavior of potentials. In Euclidean spaces, for example, 
we have \cite[Theorem 5.3 in Chapter 2]{Mi} on Riesz potentials. The following is the analog  of 
\cite[Theorem 5.3 in Chapter 2]{Mi} on potentials for sub-Laplacians. 

\begin{theorem}\label{Thm:estimate-fine-1} Let $R^{\mu, D}_\L$ be the potential of a Radon measure $\mu$ with the support in 
$D$ for the sub-Laplacian $\L$ in a homogeneous
Carnot group $\G = (\R^N, \circ, \delta_\lambda)$. Then, for any $x_0\in D$, there exists a set $E$ that is $\L$-thin at $x_0$ such that 
\begin{equation}\label{Equ:estimate-fine-1}
\lim_{x\to x_0, x\notin E} \frac {R^{\mu, D}_\L(x)}{\Gamma (x_0^{-1}\circ x)} = \mu(\{x_0\}).
\end{equation}
\end{theorem}

\begin{proof} The proof in principle 
goes with little changes from that of \cite[Theorem 5.3 in Chapter 2]{Mi}. But we want to take the opportunity and 
present the proof to show how Definition \ref{Def:thinness} is used beautifully. To deal with $d_\L$, we will need to use the reverse improved pseudo-triangle inequality, 
for example, \cite[(12.15b)]{BLU}, that is,
\begin{equation}\label{Equ:impro-pseudo-tri-r}
d_\L (y^{-1}\circ x) \geq \frac 1\beta |d_\L (x) - d_\L (y)|.
\end{equation}
For simplicity we may let $x_0=0$ and $B_\L (0, r) = B_\L(r) = \{x\in \R^N: d_\L (x)< r\}$ with respect to $d_\L$. For $x$ near $0$, we may assume $x\in \C_i$ for some $i$.  
Then we write
$$
R^{\mu, D}_\L (x) =( \int_{D\setminus B_\L(\bar\alpha^{-i_0})} + \int_{B_\L(\bar\alpha^{-i_0}) \setminus \bar\C_{i}} 
+ \int_{\bar\C_{i}})  \Gamma (y^{-1}\circ x)d\mu(y) = I + II + III.
$$
where $\bar\C_i = \C_{i-1}\cup\C_i\cup\C_{i+1}$. Hence we have, for $i_0$ to be chosen to be large but fixed,
\begin{equation}\label{Equ:I}
I \leq C \mu(D) \bar\alpha^{(Q-2)(i_0-i)} \Gamma (x) = o(1) \Gamma(x)
\end{equation}
as $x\to 0$ or equivalently $i\to\infty$. We write
$$
II = (\int_{B_\L(\bar\alpha^{-i-2})} +  \int_{B_\L(\bar\alpha^{-i_0})\setminus B_\L(\bar\alpha^{-i+1})})\Gamma (y^{-1}\circ x)d\mu(y) = II_1 + II_2
$$
where 
$$
II_1 = (\mu(\{0\}) + o(1)) \Gamma (x)  \text{ as $x\to 0$}
$$
and
$$
II_2 \leq  C(\mu(B_\L(\bar\alpha^{-i_0})) - \mu(\{0\}))\Gamma (x) = \epsilon \Gamma(x) 
$$
for any $\epsilon >0$ as $i_0$ chosen to be sufficiently large. Next, for $III$, we set
$$
E^{\lambda_i}_i = \{x\in \bar\C_i:  \int_{\bar\C_i} \Gamma(y^{-1}\circ x)d\mu(y) \geq \lambda_i \alpha^i\}
$$
for some $\lambda_i\to 0$ to be chosen. Notice that
$$
\frac { \int_{\bar\C_i} \Gamma(y^{-1}\circ x)d\mu(y)}{\Gamma(x)} \leq C \alpha^{-i} 
 \int_{\bar\C_i} \Gamma(y^{-1}\circ x)d\mu(y) \leq C \lambda_i\to 0
$$
when $x\notin E^{\lambda_i}_i$. Moreover, by Definition \ref{Def:l-capacity},
$
cap_\L^* (E^{\lambda_i}_i) \leq \frac 1{\lambda_i} \alpha^{-i}\mu(\bar\C_i).
$
Therefore, if set $E = \bigcup_i E^{\lambda_i}_i$, then
$$
\sum_n \alpha^n cap_\L^* (E\cap \C_n) \leq  \sum_n \alpha^n cap_\L^* (E^{\lambda_n}_n) \leq \sum_n \frac 1{\lambda_n}\mu (\bar\C_n).
$$
Finally, $\sum_n \alpha^n cap_\L^* (E\cap \C_n) <\infty$ for some sequence $\lambda_n\to 0$ when
$\sum_n\mu(\bar\C_n) \leq 3 \mu(D) <\infty$, thanks to \cite[(5) Page 40]{Brom} (see also  \cite{Rey}). 
Thus $E$ is $\L$-thin at the origin by Definition \ref{Def:thinness}.
\end{proof}
 
\eqref{Equ:estimate-fine-1} can be strengthened if one may choose the point $x_0$ in a given set of Hausdorff dimension $s$. 
This improvement is the key to estimate Hausdorff dimensions of singularities. Namely, we have the following extension of \cite[Theorem 1.1]{MQ-2025}. 
 
\begin{theorem}\label{Thm:theorem-7} Let $R^{\mu, D}_\L$ be the potential of a Radon measure $\mu$ supported in $D$ for the sub-Laplacian 
$\L$ in a homogeneous
Carnot group $\G = (\R^N, \circ, \delta_\lambda)$.  Let $S$ be a compact subset in $D$ whose Hausdorff dimension 
$dim_{\mathcal{H}_\L}$ is greater than $s$ for $s\in (0, Q-2)$. 
Then there is a point $x_0\in S$ and a set $E$ that is $\L$-thin at $x_0$ such that 
\begin{equation}\label{Equ:estimate-upper}
R^{\mu, D}_\L (x) \leq C(d_\L(x_0^{-1}\circ x))^{2-Q+s} 
\end{equation}
for $x$ close to $x_0$ and not in $E$.
\end{theorem} 

The proof of Theorem \ref{Thm:theorem-7} in principle goes again with little changes from that of \cite[Theorem 1.1]{MQ-2025}. 
The only issue is that we need an extension of \cite[Proposition 1.4]{Kp} to cover Carnot-Caratheodory spaces.   

\begin{lemma}\label{Lem:kpata-1.4} Let $d_\L$ be the $\L$-gauge  in \eqref{Equ:L-gauge} 
associated with a sub-Laplacian $\L$ in a homogeneous Carnot group $\G 
= (\R^N, \circ, \delta_\lambda)$. And let $\mu$ be a nonnegative Radon measure on $\R^N$ and 
$$
G^\infty_s=\{x\in \R^N: \limsup_{r\to 0} r^{-s}\mu(B_\L(x, r)) = \infty\}.
$$
Then 
$$
\mathcal{H}_\L^s (G^\infty_s) = 0
$$
where $B_\L(x, r)$ is the ball with respect to $d_\L$ and 
$\mathcal{H}_\L^s$ is the Hausdorff measure of dimension $s$ that is associated with $d_\L$ (cf. \cite[Chapter 13]{BLU}).
\end{lemma}

\begin{proof} The proof goes similar to the proof of \cite[(iii) of Proposition 1.4]{Kp}. One only needs to make sure that \cite[Lemma 2.3]{Kp} 
is applicable to Carnot-Caratheodory spaces,
which is known, for example, thanks to \cite[(1) of Theorem 2.10.18]{Fed}. 
\end{proof}

Some version of Theorem \ref{Thm:theorem-7} for Heisenberg groups appeared earlier as \cite[Lemma 1.10]{CY}.
The Radon measures $\mu$ in Theorem \ref{Thm:theorem-7}  and Lemma \ref{Lem:kpata-1.4} in this paper need not to be absolutely continuous 
with respect to the Lebesgue measure in $\R^N$. 

%%%%%%%%%%%%%%%%%%%%%%%%%%%%%%%%%%%%%%%%%%%%%%%%%%%
%%%%%%%%%%%%%%%%%%%%%%%%%%%%%%%%%%%%%%%%%%%%%%%%%%%

\section{Geodesics in Heisenberg groups and counter-examples} \label{Sect:counter-example}

In this section we construct a counter-example in Heisenberg groups to the property of thin sets 
observed in \cite[Theorem 1.1]{AH}, \cite[Lemma 8.2 in Chapter 2]{Mi}, and \cite[Theorem 1.2]{MQ-2025} 
in Euclidean cases.

%%%%%%%%%%%%%%%%%%%%%%%%%%%%%%%%%%%%%%%%%%%%%%%%%%%%

\subsection{Heisenberg groups}\label{Subsect:heisenberg-g}
The Heisenberg group, as a homogeneous Carnot group, is 
$$
\mathbb{H}^n = \{(x, y, z): x, y\in \mathbb{R}^n \text{ and } z\in \mathbb{R}\} = \R^{2n+1}
$$
with the multiplication
\begin{equation}\label{Equ:H-product}
(x, y, z)\circ (x', y', z') = (x+x', y+y', z+z' +\frac 12(x\cdot y' - x'\cdot y))
\end{equation}
and the dilation 
$
\delta_\lambda (x, y, z) = (\lambda x, \lambda y, \lambda^2 z).
$
The Jacobian generators in the Heisenberg group $\mathbb{H}^n$ are
$$
X_i = \partial_{x_i} - \frac {y_i}2 \partial_z \text{ and } Y_i = \partial_{y_i} + \frac {x_i}2 \partial_z
$$
for $i = 1, 2, \cdots, n$.
The canonical sub-Laplacian is
\begin{equation}\label{Equ:H-laplacian}
\Delta_{\mathbb{H}} = \sum_{i=1}^nX_i^2 + Y_i^2 =\Delta_X + \Delta_Y + \frac {|x|^2+|y|^2}4\partial_z^2 + (-y\cdot\partial_x+x\cdot\partial_y)\partial_z
\end{equation}
and the fundamental solution to the canonical sub-Laplacian $\Delta_\H$ is
\begin{equation}\label{Equ:folland-fundamental}
\Gamma_{\Delta_\H} (0, (x, y, z))  = [d_{\Delta_\H}(x, y, z)]^{2 - Q},
\end{equation}
where 
\begin{equation}\label{Equ:H-gauge}
d_{\Delta_\H} (x, y, z) = d_{\Delta_\H} (0, (x, y, z)) = c_n ((|x|^2+|y|^2)^2 + 16 z^2)^\frac 14
\end{equation}
is the $\Delta_{\mathbb{H}}$-gauge for some dimensional constant $c_n$ (cf.  \cite{Fo}) and 
$Q = 2n+2$ is the homogeneous dimension  of $\H^n$. The $\Delta_\H$-gauge $d_{\Delta_\H}$ satisfies the scaling property
\begin{equation}\label{Equ:scale-gauge}
d_{\Delta_\H}(\delta_\lambda (x, y, z)) = \lambda \, d_{\Delta_\H}(x, y, z) \text{ for 
$\lambda >0$}.
\end{equation}
Finally the $\Delta_\H$-gauge $d_{\Delta_\H}$ is equivalent to the Carnot-Caratheodory distance, which is the geodesic distance of the sub-Riemannian structure on $\H^n$.

%%%%%%%%%%%%%%%%%%%%%%%%%%%%%%%%%%%%%%%%%%%%%%%%

\subsection{Geodesics in Heisenberg groups} \label{Subsect:geodesics}
The description of geodesics in Heisenberg groups here follows \cite[Section 1.5 and 1.6]{RM}. 
Let us start with momentum functions
$$
P_{X_i} = p_{x_i} - \frac {y_i}2 p_z \text{ and } P_{Y_i} = p_{y_i} + \frac {x_i}2 p_z
$$
on the cotangent bundle 
$$
T^*\mathbb{H}^n = \mathbb{R}^{2n+1}\times \mathbb{R}^{2n+1} = \{(x_1, \cdots x_n, y_1, \cdots, y_n, z, p_{x_1}, \cdots p_{x_n}, p_{y_1}, 
\cdots, p_{y_n}, p_z)\}.
$$
Then the Hamiltonian for the sub-Riemannian structure associated with the Jacobian generators is
$$
H = \frac 12 \sum_{i=1}^n (P_{X_i}^2 + P_{Y_i}^2) = \frac 12\sum_{i=1}^n (p_{x_i}^2 +p_{y_i}^2 -y_ip_{x_i}p_z+x_ip_{y_i}p_z + \frac {x_i^2+y_i^2}4 p_z^2).
$$
Therefore the Hamiltonian dynamics on the cotangent bundle is
$$
\left\{
\aligned
\overset{.}{x_i}  & = P_{X_i}, \ \overset{.}{y_i}  = P_{Y_i}, \ \overset{.}{z}  = \sum_{i=1}^n (-\frac {y_i}2 P_{X_i} + \frac {x_i}2 P_{Y_i}), \\
\overset{.}{p}_{x_i}  &  = -\frac {p_z}2 P_{Y_i},  \ \overset{.}{p}_{y_i}  = \frac {p_z}2  P_{X_i}, \ \overset{.}{p}_z  = 0.
\endaligned\right.
$$
This describes geodesics of the sub-Riemannian structure on $\H^n$.

%%%%%%%%%%%%%%%%%%%%%%

\subsubsection{Normal geodesics from the origin}\label{Subsubsect:geodesics}
For geodesics from the origin with initial momentum $(p_x, p_y, 0)$,  we have
$$
x (t)  = p_x t, \ y (t)  =p_y t,\ z(t)  = 0,
$$
which is a straight ray from the origin in XY-hyperplane.  For geodesics from the origin with initial momentum $(p_x, p_y, p_z)$, where $p_z\neq 0$, 
we have
$$
\left\{\aligned
x (t) & =  \frac {p_x}{p_z} \sin p_zt + \frac {p_y}{p_z} \cos p_zt - \frac {p_y}{p_z} , \\
y (t) & = - \frac {p_x}{p_z} \cos p_zt + \frac {p_y}{p_z} \sin p_zt + \frac {p_x}{p_z},
\endaligned\right.
$$
which is a circle in XY-hyperplane centered at $(-\frac {p_y}{p_z}, \frac {p_x}{p_z})$ and with the radius $\rho= \frac{\sqrt{|p_x|^2+|p_y|^2}}{|p_z|}$
passing through the origin, and
$$
z (t)  =  \frac {|p_x|^2+|p_y|^2}{2p_z^2} (p_z t - \sin p_zt).
$$
To summarize, the geodesics from the origin are called normal geodesics and minimal, 
when they only possibly meet the $z$-axis the second time at the end. They are all horizontal lifts of circles in $XY$-hyperplane passing the origin,
including the straight lines in $XY$-hyperplane passing the origin, which are considered to be circles of an infinite radius. 
The geodesic length of the horizontal lift of a circle of a finite radius 
is the circumference of the circle, and only the straight lines in $XY$-hyperplane can go an infinite length.  

%%%%%%%%%%%%%%%%%%%%%%%

\subsubsection{The scaling property of geodesics from the origin}\label{Subsubsect:scale} As in Euclidean spaces, the dilation takes normal 
geodesics to normal geodesics and scales the length. However, in contrast to Euclidean cases, the dilation does not simply extend a normal 
geodesic from the origin, instead, it takes it to another completely distinct normal geodesic from the origin, unless the geodesic is a straight line
passing the origin inside XY-hyperplane. This distinction is significant for geometry as well as analysis on homogeneous Carnot groups.

%%%%%%%%%%%%%%%%%%%%%%%%%%%%%%%%%%%%%%%%%%%%%%%%%%%%%%

\subsection{The anisotropic nature of Heisenberg groups $\H^n$}\label{Subsect:anisotropic} 
In analysis to study potentials for sub-Laplacians on $\H^n$, we want to consider level sets of the $\Delta_\mathbb{H}$-gauge.
$$
S^{\Delta_\mathbb{H}}_r = \{(x, y, z): d_{\Delta_\H}(x, y, z) = r\}.
$$ 
Let $G^{\Delta_\mathbb{H}}_{r}$ be the collection of normal geodesics from the origin that reach to $S^{\Delta_\mathbb{H}}_r$, excluding those 
that stop at $z$-axis before they reach to $S^{\Delta_\mathbb{H}}_r$. Then, we consider the set
\begin{equation}\label{Equ:blockage}
D^{\Delta_\mathbb{H}}_{r_1, r_2} = \{(x, y, z)\in S^{\Delta_\mathbb{H}}_{r_2}: (x, y, z) \text{ on the trajectory of $\gamma\in 
G^{\Delta_\mathbb{H}}_{r_1}$}\} \text{ for $r_2< r_1$}.
\end{equation}
It is easily seen that
$
G^{\Delta_\mathbb{H}}_{r_1} \subsetneqq G^{\Delta_\mathbb{H}}_{r_2} \text{ for $r_2< r_1$}
$
and therefore
\begin{equation}\label{Equ:g-monotone}
D^{\Delta_\mathbb{H}}_{r_1', r_2} \subsetneqq D^{\Delta_\mathbb{H}}_{r_1, r_2} \text{ for $r_2 < r_1 < r_1'$}.
\end{equation}
In fact, by the summary in Section \ref{Subsubsect:geodesics}, we have
\begin{equation}\label{Equ:g-anisotropic}
D^{\Delta_\mathbb{H}}_{r, 1} \to S^{\Delta_\mathbb{H}}_1\bigcap\{z=0\} \text{ as $r\to\infty$},
\end{equation}
where $S^{\Delta_\mathbb{H}}_1\bigcap\{z=0\}$ is 
nothing but the round sphere $\mathbb{S}^{2n-1}(1/c_n)$ of radius $1/c_n$ in the XY-hyperplane by 
\eqref{Equ:H-gauge}. Finally we want to mention that
$$
D^{\Delta_\mathbb{H}}_{\lambda r_1, \lambda r_2} = \delta_\lambda (D^{\Delta_\mathbb{H}}_{r_1, r_2}) \text{ for $\lambda>0$},
$$
due to the scaling property of both $\Delta_\H$ gauge \eqref{Equ:scale-gauge} and normal geodesics (cf. Section \ref{Subsubsect:scale}).

%%%%%%%%%%%%%%%%%%%%%%%%%%%%%%%%%%%%%%%%%

\subsection{Geodesically unavoidable $\Delta_\H$-thin sets}\label{Subsect:counter-ex}  
Let $cap_{\Delta_\H}^*$ be the exterior capacity defined in Definition \ref{Def:l-capacity} associated with the canonical sub-Laplacian $\Delta_\H$ in Heisenberg group $\H^n$. 
First, in the light of the scaling property \eqref{Equ:Gamma-scale} of fundamental solution $\Gamma_{\Delta_\H}$,  we have 
\begin{equation}\label{Equ:scaling}
cap_{\Delta_\H}^*(\delta_\lambda(E)) = \lambda^{Q-2} cap_{\Delta_\H}^* (E).
\end{equation}

\begin{lemma} Let $\mathbb{S}^{2n-1}(1/c_n) = S^{\Delta_\H}_1 \bigcap\{z=0\}$. Then
$$
 cap_{\Delta_\H}^* (S^{\Delta_\H}_1 \cap \{z=0\}) = 0.
$$
\end{lemma}

\begin{proof} For simplicity, we want to show $cap_{\Delta_\H}^*(\mathbb{S}^{2n-1}) = 0$ instead
for the unit round sphere $\mathbb{S}^{2n-1}$ in $XY$-hyperplane of $\H^n$ due to \eqref{Equ:scaling}. 
According to \cite[Lemma 11.5.11]{BLU}, it suffices show that $\mathbb{S}^{2n-1}$ is $\Delta_\H$-polar. Let
$\mu_0(A) = m_0 (A\bigcap \mathbb{S}^{2n-1})$
be the Radon measure on $\H^n$, where $m_0$ is the standard Lebesgue measure on $\mathbb{S}^{2n-1}$. Then we calculate 
the potential
\begin{equation}\label{Equ:potential-polar}
R^{\mu_0, \R^{2n+1}}_{\Delta_\H} (x, y, 0) = \int_{\mathbb{S}^{2n-1}} \Gamma_{\Delta_\H} ((x, y, 0)^{-1}\circ u) dm_0(u) 
= \infty 
\end{equation}
for $(x, y, 0) \in \mathbb{S}^{2n-1}\subset\H^n$.
In fact, due to the symmetry, one only needs to calculate the potential $R^{\mu_0, \R^{2n+1}}_{\Delta_\H}$
at $x_0=(1, 0, \cdots, 0)$ and $y_0 = (0, 0, \cdots, 0)$. Notice that
$$
\Gamma_{\Delta_\H} ((x_0, y_0, 0)^{-1}\circ u) = c_n^{-2n} 2^{-n}((1 - u_{x_1})^2 + u_{y_1}^2)^{-\frac n2}
$$
for $u=(u_x, u_y, 0) \in \mathbb{S}^{2n-1}\subset \H^n$, by \eqref{Equ:H-product}, \eqref{Equ:folland-fundamental}, and \eqref{Equ:H-gauge}.
To complete \eqref{Equ:potential-polar}, one uses the spherical coordinates on $\mathbb{S}^{2n-1}$, in which 
\begin{itemize}
\item $u_{x_1} = \cos\theta$, $u_{y_1} = \sin\theta$, where $\theta\in[0, 2\pi]$, and
$
R^{\mu_0, \R^{2n+1}}_{\Delta_\H} (x_0, y_0, 0)\geq c \int_0^\delta \theta^{-1}d\theta = \infty 
$
for some $\delta >0$ and $c>0$, when $n=1$;
\item $u_{x_1} = \cos\theta_1$, $u_{y_1} = \sin\theta_1\sin\theta_2$,  where $\theta_1\in[0, \pi]$ and $\theta_2 \in [-\frac \pi 2, \frac \pi 2]$, and
$$
R^{\mu_0, \R^{2n+1}}_{\Delta_\H}  (x_0, y_0, 0)\
 \geq c\int_0^\delta \theta_1^{n-2} \int_0^{\theta_1} (\theta_1^2 + \theta_2^2)^{-\frac n2} \, d\theta_2\, d\theta_1
\geq  2^{-\frac n2}c\int_0^\delta \theta_1^{-1}\, d\theta_1 = \infty
$$
\end{itemize} 
for some $\delta > 0$ and $c>0$, when $n\geq 2$.
\end{proof}

Consequently, by (11.19) in  \cite[Lemma 11.5.10]{BLU} as well as \eqref{Equ:g-monotone} and \eqref{Equ:g-anisotropic}, we have

\begin{lemma} 
$
\lim_{r\to\infty} cap_{\Delta_\H}^* (D^{\Delta_\mathbb{H}}_{r, 1}) = 0.
$
\end{lemma}

Therefore we may select a sequence $k_i\to\infty$ such that
$\sum_{i=1}^\infty cap_{\Delta_\H}^* (D^{\Delta_\mathbb{H}}_{\bar\alpha^{k_i}, 1}) < \infty$ and set
\begin{equation}\label{Equ:set-E}
E =  \bigcup_{i=1}^\infty  D^{\Delta_\mathbb{H}}_{\bar\alpha^{-k_i}, \bar\alpha^{-2k_i}}.
\end{equation}

\begin{proposition} \label{Prop:thin-example} For the set $E$ given in \eqref{Equ:set-E},
$\sum_{n=1}^\infty \alpha^{n} cap_{\Delta_\H}^* (E\bigcap\C_n) < \infty.$
In other words, the set $E$ in \eqref{Equ:set-E} is $\Delta_\H$-thin at the origin according to Definition \ref{Def:thinness}.
\end{proposition}

\begin{proof} First, from \eqref{Equ:set-E}, 
$$
\sum_{n=1}^\infty \alpha^n  cap_{\Delta_\H}^* (E\bigcap\C_n)  = \sum_{i=1}^\infty \alpha^{2k_i} cap_{\Delta_\H}^* (E\bigcap\C_{2k_i}) 
 = \sum_{i=1}^\infty \alpha^{2k_i} cap_{\Delta_\H}^* ( D^{\Delta_\mathbb{H}}_{\bar\alpha^{-k_i}, \bar\alpha^{-2k_i}}).
$$
Then, by the scaling property \eqref{Equ:scaling} and \eqref{Equ:alpha-bar},
$$
 \sum_{i=1}^\infty \alpha^{2k_i} cap_{\Delta_\H}^* ( D^{\Delta_\mathbb{H}}_{\bar\alpha^{-k_i}, \bar\alpha^{-2k_i}}) 
 \leq 3 \sum_{i=1}^\infty cap_{\Delta_\H}^* (D^{\Delta_\mathbb{H}}_{\bar\alpha^{k_i}, 1}) < \infty.
$$
\end{proof}

Consequently,

\begin{theorem}\label{Thm:E-counter} In Heisenberg group $\mathbb{H}^n$, the set $E$ in \eqref{Equ:set-E}
is $\Delta_\H$-thin at the origin according to Definition \ref{Def:thinness}. However, any geodesic from the origin intersects with $E$.
\end{theorem}
\begin{proof} This is simply because that a geodesic is always normal initially and that a normal geodesic always 
intersects with $E$. In fact, due to the definition \eqref{Equ:blockage}, any normal geodesic that reaches to the level set of $S^{\Delta_\H}_{\bar\alpha^{-k_i}}$ intersects 
with $D_{\bar\alpha^{-k_i}, \bar\alpha^{-2k_i}}$ as $k_i\to\infty$.
\end{proof}

%%%%%%%%%%%%%%%%%%%%%%%%%%%%%%%%%%%%%%%%%%%%%
%%%%%%%%%%%%%%%%%%%%%%%%%%%%%%%%%%%%%%%%%%%%%

\section{On thin sets in polarizable groups and geometric applications}\label{Sect:escape}

It is worth pointing out that, for a version of the existence of rays escaping from a thin set 
in Carnot groups to utilize the geometric completeness, one may use horizontal curves which need not be geodesics. 
Also a measurement for a set to be of small size in $\textup{S}\setminus\mathcal{Z}$ 
may simply be a measure instead of the capacity.  
With those in mind, we will employ the foliation of horizontal radial curves for the polar coordinates in polarizable
Carnot groups constructed in \cite{BT} to establish Theorem \ref{Thm:capacity-estimate}.  
$\G$ is assumed to be polarizable (cf. \cite[Definition 2.12]{BT}) throughout this section. 

%%%%%%%%%%%%%%%%%%%%%%%%%%%%%%%%%%%%%%%%%%%%%

\subsection{Polarizable groups and polar coordinates}\label{Subsect:polar-coordinate} In this subsection we follow \cite{BT} to introduce polar
coordinates using the foliation of horizontal radial curves in polarizable Carnot groups. 
To introduce flow lines as the horizontal radial curves, we let 
$$
\textup{S} = \{x\in \R^N: d_{\Delta_\G}(x) = 1\} \text{ and } \mathcal{Z} = \{0\}\bigcup \{x\in\R^N: \nabla_0 d_{\Delta_\G} (x) = 0\},
$$ 
where $\nabla_0 d_{\Delta_\G}$ is the horizontal gradient of the $\Delta_\G$-gauge $d_{\Delta_\G}$. And in $\R^N\setminus\mathcal{Z}$, one considers the 
flow lines $\phi(s, g):  \R^N\setminus\mathcal{Z}\to  \R^N\setminus\mathcal{Z}$ that satisfy
\begin{equation} \label{Equ:flow-BT}
\left\{ \aligned \frac {\partial}{\partial s} \phi(s, g) & = \frac {d_{\Delta_\G}(\phi(s, g))}s \cdot 
\frac {\nabla_0 d_{\Delta_\G} (\phi(s, g))} {\|\nabla_0 d_{\Delta_\G}(\phi(s, g))\|_0^2} \\
\phi(1, g)  & = g \endaligned\right.,
\end{equation}
where $\|\cdot\|_0$ is the norm on the horizontal space with respect to which the Jacobian generators form the orthonormal basis.

\begin{lemma}\label{Lem:polar-coordinate} (\cite[Lemma 3.3]{BT}) Let $\G= (\R^N, \circ, \delta_\lambda)$ be polarizable. 
Then the flow lines $\phi(s, g)$ satisfies
\begin{itemize}
\item $d_{\Delta_\G}(\phi(s, g)) = s d_{\Delta_\G}(g)$ for $s>0$ and $g\in \R^N\setminus\mathcal{Z}$;
\item the velocity  $\frac {\partial}{\partial s} \phi(s, g)$ is horizontal and the speed 
$
\|\frac\partial{\partial s}\phi (s, g)\|_0 = \frac {d_{\Delta_\G}(g)}{\|\nabla_0d_{\Delta_\G}(g)\|_0}
$ 
is independent of $s$.
%\item $\det D \phi(s, g) = s^Q$ for $s>0$ and $g\in \R^N\setminus\mathcal{Z}$.
\end{itemize}
\end{lemma}

\begin{lemma}\label{Lem:polar-integral} (\cite[Theorem 3.7]{BT}) Let $\G = (\R^N, \circ, \delta_\lambda)$ be polarizable.
There exists a unique Radon measure $\sigma$ on $\textup{S}\setminus\mathcal{Z}$ such that for all $u\in L^1(\R^N)$
\begin{equation}\label{Equ:polar-integral}
\int_{\R^N} u(x) dx = \int_0^\infty \int_{\textup{S}\setminus\mathcal{Z}} u(\phi(s, v)) s^{Q-1} d\sigma(v) ds.
\end{equation}
\end{lemma}

Thanks to \cite[Proposition 5.6]{BT}, H-type groups (cf. \cite[Definition 5.1]{BT}) are polarizable. And, 
according to \cite{CDKR-91}, H-type groups include those that are the nilpotent part of Iwasawa decomposition 
of rank one groups: $SU(n+1, 1)$, $Sp(n+1, 1)$ and $\mathbb{F}_{4(-20)}$, which are the relevant groups for 
our geometric applications.

%%%%%%%%%%%%%%%%%%%%%%%%%%%%%%%%%%%%%%%%%%%%%%%

\subsection{The flow lines $\phi(s, \cdot)$ as horizontal radial curves} \label{Subsect:key-lemma}
To fix the upper bound for the speed of flow lines $\phi(s, \cdot)$ in the light of
Lemma \ref{Lem:polar-coordinate}, let us introduce
\begin{equation}\label{Equ:a-S}
\textup{S}_a = \{g\in \textup{S}: \|\nabla_0d_{\Delta_\G}(g)\|_0 > a\} \text{ and } \G_a = \{\phi(s, g): g\in \textup{S}_a \text{ and } s\in (0, \infty)\}.
\end{equation}

Then

\begin{lemma} Let $\G$ be a polarizable group and $g\in \textup{S}_a$ for $a>0$. Then
\begin{equation}\label{Equ:speed}
d_{\Delta_\G}(\phi(s, g)^{-1}\circ\phi(r, g)) \leq C_a |r-s|
\end{equation}
for some constant $C_a$, which depends only on the number $a$.
\end{lemma}

\begin{proof} It suffices to prove \eqref{Equ:speed} for the Carnot-Caratheodory norm $d_c$ instead of the gauge $d_{\Delta_\G}$. Because they are equivalent on any Carnot group.
To finish the proof, one may use the speed upper bound of the horizontal radial curve $\phi(t, g)$ connecting $\phi(r, g)$ and $\phi(s, g)$. 
\end{proof}

Let
$$
\mathcal{P} (x) = \mathcal{P} (\phi(s, g)) = g: \R^N\setminus\mathcal{Z} \to S \setminus\mathcal{Z}
$$
be the projection based on the polar coordinates by horizontal radial curves $x = \phi(s, g)$ for $g\in \textup{S}\setminus\mathcal{Z}$.
The following is the key lemma of this note that facilitates the existence of horizontal radial curves escaping from a thin set. 
From now on we will always work with $\C_n$ that are centered at the origin.

\begin{lemma}\label{Lem:cap-by-area} Suppose that $\mathbb{G} = (\R^N, \circ, \delta_\lambda)$ is polarizable. 
Let $F$ be an open set and $a>0$ be appropriately small. Then
\begin{equation}\label{Equ:cap-by-area}
cap_{\Delta_\G}^* (F\cap\C_n\cap\G_a) \geq C \bar\alpha^{-(Q-2)n}\sigma (\mathcal{P}(F\cap\C_n\cap\G_a))^\frac {Q-2}{Q-1},
\end{equation}
where the measure $\sigma$ is the one in Lemma \ref{Lem:polar-integral} and $C$ depends on $a$.
\end{lemma}

Before we prove this key lemma, we first use it to prove the existence of horizontal radial curves escaping from a thin set
in polarizable Carnot groups.

\begin{theorem} \label{Thm:capacity-estimate} Let $\G = (\R^N, \circ, \delta_\lambda)$ be polarizable. And let $E$ be a set that is thin at the origin in 
$\G$. Then there always exists a flow line $\phi(s, g)$ given by \eqref{Equ:flow-BT} that does not intersect with $E$, at least in a small neighborhood of the
origin.
\end{theorem}

\begin{proof}
According to Definition \ref{Def:thinness}, we consider, for a thin set $E$, 
\begin{equation}
\aligned
\infty  & > \sum_{n=1}^\infty \bar\alpha^{n(Q-2)} cap_{\Delta_\G}^* (E\cap \C_n) \geq  \sum_{n=1}^\infty \bar\alpha^{n(Q-2)} cap_{\Delta_\G}^* (E\cap \C_n\cap\G_a) \\
& \geq \frac 12   \sum_{n: \, cap_{\Delta_\G}^*(E\cap\C_n\cap\G_a) > 0} 
\bar\alpha^{n(Q-2)} cap_{\Delta_\G}^*(F\cap \C_n\cap\G_a) \\
& \geq C \sum_{n: \, cap_{\Delta_\G}^*(E\cap\C_n\cap\G_a) > 0} (\sigma (\mathcal{P}(F\cap \C_n\cap\G_a)))^\frac {Q-2}{Q-1}\\
&  \geq C \sum_{n: \, cap_{\Delta_\G}^*(E\cap\C_n\cap\G_a) > 0} \sigma (\mathcal{P}(F\cap \C_n\cap\G_a)) \\
& \geq C \sigma(\bigcup_{n: \, cap_{\Delta_\G}^*(E\cap\C_n\cap\G_a) > 0} \mathcal{P}(F\cap\C_n\cap\G_a))
\endaligned
\end{equation}
and
\begin{equation}
\aligned
1 & = \sum_{n=1}^\infty 2^{-n} 
 \geq \sum_{n: \, cap_{\Delta_\G}^*(E\cap\C_n\cap\G_a) = 0} 
 \bar\alpha^{n(Q-2)} cap_{\Delta_\G}^*(F\cap \C_n\cap\G_a)\\
& \geq \sigma(\bigcup_{n: \, cap_{\Delta_\G}^*(E\cap\C_n\cap\G_a) = 0} \mathcal{P}(F\cap\C_n\cap\G_a)) 
\endaligned
\end{equation}
for som appropriate open set $F\supset E$ due to the definition of the exterior capacity (cf. Definition \ref{Def:l-capacity}). This implies $\textup{S}_a\setminus \bigcup_{n=k}^\infty \mathcal{P}(F\cap\C_n\cap\G_a)$ is not empty when 
$k$ is sufficiently large. Because
\begin{equation}
\sigma({\bigcap_{k=1}^\infty} \bigcup_{n=k}^\infty \mathcal{P}(F\cap \C_n\cap\G_a)) = 0 < \sigma(\textup{S}_a),
\end{equation}
which is like \cite[Lemma 8.2 in Chapter 2]{Mi} and suffices for our purpose. $\sigma(\textup{S}_a)$ is positive when $a$ is sufficiently small 
(see \cite[Proposition 2.18 (i), Lemma 3.3, and Theorem 3.7]{BT}). 
\end{proof}

We now prove Lemma \ref{Lem:cap-by-area} in the following subsections. The idea is to construct a competitive 
Radon measure supported in $F\cap\C_0\cap\G_a$ for open sets $F$. Then Lemma \ref{Lem:cap-by-area} in general follows in the light of the 
scaling property of the homogeneously induced measure on $\phi(s, \textup{S}\setminus\mathcal{Z})$:
\begin{equation}\label{Equ:s-measure}
\sigma_s(\phi(s, A)) = s^{Q-1} \sigma (A)
\end{equation}

%%%%%%%%%%%%%%%%%%%%%%%%%%

\subsubsection{A competitive Radon measure}\label{Subsubsect:competitive-measure} For a set $F\cap\mathcal{C}_0\cap\G_a$ and 
a Radon measure $\mu$ supported in $F\cap\C_0\cap\G_a$ is said to be competitive if it satisfies
\begin{enumerate}
\item For a positive constant $A$
\begin{equation}\label{Equ:lower-b}
\sigma(\mathcal{P}(F\cap\C_0\cap\G_a)) \geq \mu (F\cap\C_0\cap\G_a) \geq A\sigma(\mathcal{P}(F\cap\C_0\cap\G_a)).
\end{equation}
\item For $\delta_0>0$ there exists a constant $B$  
\begin{equation}\label{Equ:upper-b}
\mu (B_{\Delta_\G}(p, \rho)\cap F\cap\C_0\cap\G_a) \leq B \rho^{Q-1} 
\end{equation}
 when $d_{\Delta_\G}(p) = r \geq \delta_0$,  where $B_{\Delta_\G}(p, \rho) = \{g\in G: d_{\Delta_\G}(p^{-1}\circ g) <\rho\}$.
 \end{enumerate}
 
 %%%%%%%%%%%%%%%%%%%%%%%%%%%%
 
\subsubsection{A competitive Radon measure for spherical wedge shells}\label{Subsubsect:measure-shell} 
It is easy to have a competitive measure for a spherical wedge shell with respect to the flow lines $\phi(s, \cdot)$. 
Here a spherical wedge shell with respect to the flow lines $\phi(s, \cdot)$ is
\begin{equation}\label{Equ:wedge-shell}
W^A_{r_1, r_2} = \{ \phi(s, g)\in \G: g\in A\subset \textup{S}_a \text{ and } s\in (r_1, r_2]\} = \bigcup_{s\in (r_1, r_2]}\phi(s, A),
\end{equation}
where $\phi(s, A) = \{\phi(s, g): g\in A\subset \textup{S}_a\}$ is a cap in $\G_a$ for a measurable set $A$ and $s>0$. Clearly we have
$\mathcal{P}(W^A_{r_1, r_2}) = A.$ A candidate for competitive measures is simply the restriction 
$$
\mu_{r, A} (E) =  r^{1-Q} \sigma_r (E \cap \phi(r, A)) = \sigma(\mathcal{P}(E\cap \phi(r, A))
$$
for a fixed $r\in (r_1, r_2]\subset (\bar\alpha^{-1}, 1]$, which clearly is a Radon measure supported on the cap $\phi(r, A)$ in 
$W^A_{r_1, r_2}$ and satisfies \eqref{Equ:lower-b}. 
In order to verify \eqref{Equ:upper-b} for $\mu_{r, A}$, we recall the improved pseudo-triangle inequality \cite[(12.15a)]{BLU}
\begin{equation}\label{Equ:improved-tri}
d_{\Delta_\G}(u^{-1}\circ w) = d_{\Delta_\G}((u^{-1}\circ v)\circ (v^{-1}\circ w)) \leq \beta d_{\Delta_\G}(u^{-1}\circ v) + d_{\Delta_\G}(v^{-1}\circ w)
\end{equation}
for some $\beta \geq 1$. Then we observe

\begin{lemma}\label{Lem:upper-estimate}  
Let $\Sigma_{p, \rho}= \mathcal{P}(B_{\Delta_\G}(p, \rho)\cap\G_a)\subset \textup{S}_a$ for 
$d_{\Delta_\G}(p) =r >0$.  Then
\begin{equation}\label{Equ:inscribed}
\phi(t, \Sigma_{p, \rho}) \subset B_{\Delta_\G}(p, (1+2C_a)\beta\rho) \text{ for $t\in (r, r+ \rho)$}.
\end{equation}
Consequently we have
\begin{equation}\label{Equ:upper-b-ball}
\sigma(\Sigma_{p, \rho}) \leq \frac {\omega_\G((1+2C_a)\beta)^Q}{r^{Q-1}} \rho^{Q-1},
\end{equation}
where $\omega_\G$ is the Euclidean volume of $B_{\Delta_\G}(0, 1)$ in $\G$ and $C_a$ is from \eqref{Equ:speed}. 
\end{lemma}
\begin{proof} From \eqref{Equ:improved-tri}, for any $\phi(s, g)\in B_{\Delta_\G}(p, \rho)\cap\G_a$ and $g\in \Sigma_{p, \rho}\subset\textup{S}_a$, one has
$$
d_{\Delta_\G}(p^{-1}\circ \phi(r, g))  \leq \beta d_{\Delta_\G}(p^{-1}\circ\phi(s, g)) + d_{\Delta_\G}(\phi(s, g)^{-1}\circ\phi(r, g)) \leq  \beta \rho + C_a|r-s|,
$$
where $C_a$ is from \eqref{Equ:speed}. On the other hand, from \eqref{Equ:impro-pseudo-tri-r}, one also has
$$
d_{\Delta_\G}(p^{-1}\circ \phi(s, g))\geq \frac 1\beta |r - s|
$$
which implies $|r-s| \leq \beta \rho$ if $\phi(s, g)\in B_{\Delta_\G}(p, \rho)$. Therefore
$$
d_{\Delta_\G}(p^{-1}\circ \phi(r, g)) \leq (1 + C_a)\beta \rho \text{ for all $g\in \Sigma_{p, \rho}$.}
$$
Applying \eqref{Equ:improved-tri} again, one gets, for $g\in \Sigma_{p, \rho}$,
$$
d_{\Delta_\G}(p^{-1}\circ\phi(t, g))\leq \beta C_a |t-r| + (1+C_a)\beta \rho < (1+2C_a)\beta \rho \text{ when $t\in(r,  r+\rho)$.}
$$ 
Therefore \eqref{Equ:inscribed} is verified. 
To complete the proof, we invoke Lemma \ref{Lem:polar-integral} and the fact that the Euclidean measure is the Haar measure for $\G$. Therefore
$$
\omega_\G ((1+2C_a)\beta)^Q \rho^Q 
 = \int_{B_{\Delta_\G}(p, \, (1+2C_a)\beta\rho)} dx \geq \int_{r}^{r+\rho}\int_{\Sigma_{p, \rho}} s^{Q-1}d\sigma ds  
 \geq r^{Q-1} \rho \sigma(\Sigma_{p, \rho}).
$$
\end{proof}

Consequently \eqref{Equ:upper-b} for $\mu_{r, A}$ holds from Lemme \ref{Lem:upper-estimate}.

\begin{corollary}\label{Cor:upper-estimate} Let $d_{\Delta_\G}(p) = r\geq \delta_0>0$ and $W^A_{r_1, r_2}\subset\C_0\cap\G_a$.
Then
$$
\mu_{r, A} (B_{\Delta_\G}(p, \rho)\cap W^A_{r_1, r_2}) \leq \sigma(\mathcal{P}(B_{\Delta_\G}(p, \rho) \cap\G_a))\leq C \rho^{Q-1}.
$$
\end{corollary}
\begin{proof} Under the assumption that $d_{\Delta_\G}(p) = r \geq \delta_0>0$, \eqref{Equ:upper-b-ball} implies
\begin{equation}\label{Equ:estimate-measure}
\mu_{r, A}(W^A_{r_1, r_2}\cap B_{\Delta_\G}(p, \rho)) \leq \sigma(A\cap\Sigma_{p, \rho}) 
\leq \sigma(\Sigma_{p, \rho})\leq C\rho^{Q-1}
\end{equation}
for some constant $C$ determined by $\delta_0$ and $a$.
\end{proof}

%%%%%%%%%%%%%%%%%%%%%%%%%%%%%%

\subsubsection{A competitive measure for sets $F\cap\C_0\cap\G_a$}\label{Subsubsect:measure-open}

To construct a competitive measure for avset $F\cap\C_0\cap\G_a$, where $F$ is open, we want to pick up 
adequately many spherical wedge shells inside $F\cap\C_0\cap\G_a$ such that the union of their
projections approaches the projection of $A$ in measure in $\textup{S}\setminus\mathcal{Z}$. The strategy is to use 
the grids of thinner and thinner spherical shells to divide $A$. 

For natural numbers $k= 1, 2, \cdots$ and $l= 0, 1, \cdots, 2^k-1$, let $s^k_l = \bar\alpha^{-1+2^{-k}l}$ and let 
\begin{equation}\label{Equ:s-shells}
\mathcal{S}_l^k = \{\phi(s, g): s\in (s^k_l, s^k_{l+1}] \text{ and } g\in\textup{S}\setminus\mathcal{Z}\} 
= \{x\in \R^N\setminus\mathcal{Z}: d_{\Delta_\G}(x)\in (s^k_l, s^k_{l+1}]\}
\end{equation}
be the spherical shells with respect to the gauge $d_{\Delta_\G}$. 
For an open set $A\subset \C_0\cap\G_a$, we denote
$$
A^k_l= A\cap \mathcal{S}^k_l \text{ and } A_t = \{\phi(s, g)\in A: s=t \text{ and } g\in \textup{S}_a\}.
$$
Then we consider
$$
\underline{A}_l^k  =  \{\phi(s, g): g\in\bigcap_{t\in (s_l^k, s_{l+1}^k]} \mathcal{P}(A_t) \text{ and } s\in (s_l^k, s_{l+1}^k]\},
$$
which is the largest spherical wedge shell inside $A_l^k$ of the radius between $s^k_l$ and $s^k_{l+1}$ at $k$th level 
with respect to the flow lines $\phi(s, \cdot)$. Hence
$
\underline{A}^k = \bigcup_{l=0}^{2^k-1} \underline{A}_l^k \subset A,
$
where $\underline{A}^k$ is a union of spherical wedge shells inside $A$ at $k$th level and monotone increasing as $k$ increases. To avoid overlaps or multiplicities, 
for each $k$, we let
$$
\underline{B}^k_l = \{\phi(s, g): g\in (\bigcap_{t\in (s^k_l, s^k_{l+1}]}\mathcal{P}(A_t))\setminus \bigcup_{m=0}^{l-1} \mathcal{P}(\underline{A}^k_m) \text{ and } s\in (s^k_l, s^k_{l+1}]\}
\subset \underline{A}^k_l
$$
inductively from $l=0$ to $l=2^k-1$ ($\underline{B}^k_l$ may be empty). And we let
$
\underline{B}^k = \bigcup_{l=0}^{2^k-1} \underline{B}^k_l \subset \underline{A}^k\subset A.
$
This makes 
$$
\bigcup_{l=0}^{2^k-1}\mathcal{P}(\underline{A}^k_l) =  \mathcal{P}(\underline{A}^k) = \mathcal{P}(\underline{B}^k) = \bigcup_{l=0}^{2^k-1} \mathcal{P}(\underline{B}^k_l)
$$ 
is a union of disjoint subsets in $\textup{S}\setminus\mathcal{Z}$ and $\{\underline{B}^k_l\}_{l=0}^{2^k-1}$ is a collection of disjoint spherical wedge shells in $A$ for each $k$. 

\begin{lemma}\label{Lem:filling-by-shell} Suppose that $A = F\cap\C_0\cap\G_a$, where $F$ is open.   
Then $\{\underline{A}^k\}$ is an increasing sequence of unions of spherical wedge shells inside $A$ as $k$ increases and 
\begin{equation}\label{Equ:direct-limit}
\bigcup_{k=1}^\infty \underline{A}^k = A.
\end{equation}
Hence
$
\sigma(\mathcal{P}(\underline{A}^k))
\to \sigma(\mathcal{P}(A)) \text{ as $k\to\infty$},
$
where $\sigma$ is the Radon measure on $\textup{S}\setminus\mathcal{Z}$ in \eqref{Equ:polar-integral}. Meanwhile, for each $k$, we have
\begin{equation}\label{Equ:disjoint}
\bigcup_{l=0}^{2^k-1}\mathcal{P}(\underline{B}^k_l) = \mathcal{P}(\underline{B}^k)  = \mathcal{P}(\underline{A}^k).
\end{equation}
Thus $\sum_{l=0}^{2^k-1} \sigma(\mathcal{P}(\underline{B}^k_l)) = \sigma(\mathcal{P}(\underline{B}^k)) = \sigma(\mathcal{P}(\underline{A}^k))\to 
\sigma(\mathcal{P}(A)) \text{ as $k\to\infty$.}$
\end{lemma}

\begin{proof} It suffices to see \eqref{Equ:direct-limit} holds because that $A$ is open in $\C_0\cap\G_a$ with respect to the induced topology on $\C_0\cap\G_a$.
\end{proof}

For each spherical wedge shell $\underline{B}^k_l$, we pick the Radon measure $\mu_{r^k_l, \mathcal{P}(\underline{B}^k_l)}$ for some $r^k_l\in (s^k_l, s^k_{l+1}]$ as in Section \ref{Subsubsect:measure-shell}.  Let
\begin{equation}\label{Equ:measure-k}
\mu_k = \sum_{l=0}^{2^k-1}\mu_{r^k_l, \mathcal{P}(\underline{B}^k_l)}.
\end{equation}
Then 
\begin{equation}\label{Equ:mu-k-good}
\mu_k(A) = \sum_{l=0}^{2^k-1} \sigma(\mathcal{P}(\underline{B}^k_l)) =\sigma(\mathcal{P}(\underline{B}^k)) = \sigma(\mathcal{P}(\underline{A}^k)) \to \sigma(\mathcal{P}(A)) \text{ as $k\to\infty$.}
\end{equation}
To summarize we have 

\begin{lemma}\label{Lem:measure-open} Suppose that $A = F\cap\C_0\cap\G_a$, where $F$ is open.
Then the Radon measure $\mu_k$ given by \eqref{Equ:measure-k} supported in $A$ satisfies both \eqref{Equ:lower-b} 
and \eqref{Equ:upper-b} when $k$ is sufficiently large.
\end{lemma}

\begin{proof} \eqref{Equ:lower-b} holds for $\mu_k$ for sufficiently large $k$ by \eqref{Equ:mu-k-good}. And \eqref{Equ:upper-b} holds for $\mu_k$ because of Lemma 
\ref{Lem:upper-estimate} and the last inequality in \eqref{Equ:estimate-measure} in the proof of Corollary \ref{Cor:upper-estimate}. Because $W^A_{r_1, r_2}$ and $A$ are not
really essential in \eqref{Equ:estimate-measure}.
\end{proof}

%%%%%%%%%%%%%%%%%%%%%%%%%%%%%%%%

\subsubsection{The estimate of capacity  - the proof of Lemma \ref{Lem:cap-by-area}}\label{Subsubsect:proof} 
To estimate the capacity $cap_{\Delta_\G}^*(A)$, where $A=F\cap\C_0\cap\G_a$ and $F$ is open, 
we  consider the potential of the measure $\mu_k$, for sufficiently large $k$, constructed in Section \ref{Subsubsect:measure-open}. To invoke Definition \ref{Def:l-capacity}, let $K$ be
any compact subset of $A$ such that 
\begin{equation}\label{Equ:mu-k-lower-b}
\sigma(\textup{S}\setminus\mathcal{Z})\geq \mu_k(A) \geq \mu_k(K) \geq \frac 12\mu_k(A)\geq C\sigma(\mathcal{P}(A)).
\end{equation}
Now, let us calculate the potential for the Radon measure $\mu_k$ restricted to $K$. We claim
\begin{equation}\label{Equ:claim}
R^{\mu_k, K}_{\Delta_\G} (p) \leq C \mu_k(K)^\frac 1{Q-1}
\end{equation}
for constant $C$ and all $p\in \R^N$. For $p\in \{p\in\R^N: \inf_{q\in K}d_{\Delta_\G}(p^{-1}\circ q) \geq \mu_k(K)^\frac 1{Q-1}\}$,
it is easily calculated that 
$$
R^{\mu_k, K}_{\Delta_\G} (p)  = \int_K \frac 1{d_{\Delta_\G}(p^{-1}\circ q)^{Q-2}}d\mu_k(q) \leq \mu_k(K)^\frac 1{Q-1}. 
$$
Otherwise, we use \eqref{Equ:upper-b} to calculate
$$
\aligned
& R^{\mu_k, K}_{\Delta_\G} (p)  = \int_K \frac 1{d_{\Delta_\G}(p^{-1}\circ q)^{Q-2}}d\mu_k(q) = (Q-2) 
 \int_0^\infty \frac {\mu_k(B_t(p)\cap K)}{t^{Q-2}} \frac {dt}{t} \\
& \leq (Q-2)(B \int_0^{\mu_k(K)^\frac 1{Q-1}} dt + \mu_k(K) \int_{\mu_k(K)^\frac 1{Q-1}}^\infty t^{1-Q} dt) \\
& \leq (B(Q-2) + 1)\mu_k(K)^\frac 1{Q-1},
\endaligned
$$
provided that 
\begin{equation}\label{Equ:good-case}
\mu_k(K)^\frac 1{Q-1}  < \frac 12\bar\alpha^{-1}.
\end{equation}
Because, for some $q\in K\subset \C_0$, 
$d_{\Delta_\G}(p^{-1}\circ q)<\mu_k(K)^\frac 1{Q-1}$, and therefore, by \eqref{Equ:impro-pseudo-tri-r},
$$
d_{\Delta_\G}(p) \geq \frac 1\beta |d_{\Delta_\G}(q) - d_{\Delta_\G}(p^{-1}\circ q)| \geq \frac 1{2\beta} \bar\alpha^{-1} >0.
$$
This makes \eqref{Equ:upper-b} is applicable. If \eqref{Equ:good-case} is reversed, then we repeat the above calculations after
substituting $\mu_k(K)^\frac 1{Q-1}$ with $\frac 12\bar\alpha^{-1}$ and arrive at
$$
R^{\mu_k, K}_{\Delta_\G}(p) \leq C\mu_k(K)^\frac 1{Q-1} + C \bar\alpha^{Q-2} \mu_k(K)
$$
which implies \eqref{Equ:claim} since $\mu_k(K)\leq \mu_k(A) \leq \sigma(\textup{S}\setminus\mathcal{Z})<\infty$. 

Thus, for the measure $\tilde \mu_k$ that is the measure $\mu_k$ rescaled by the right hand side of \eqref{Equ:claim}, we have
$$
R^{\tilde\mu_k, K}_{\Delta_\G} (p) \leq 1 \text{ for all $p\in \R^N$}
$$
and 
$$
cap_{\Delta_\G}(K) \geq \tilde\mu_k(K) = C \mu_k(K)^\frac {Q-2}{Q-1}\geq C\sigma(\mathcal{P}(A))^\frac {Q-2}{Q-1}
$$
in the light of Definition \ref{Def:l-capacity} and \eqref{Equ:mu-k-lower-b}. 
The proof of  Lemma \ref{Lem:cap-by-area} is complete.

%%%%%%%%%%%%%%%%%%%%%%%%%%%%%%%%%%%%%%%%%%%%%%%%%%%%%

\subsection{Geometric applications}\label{Subsect:geometric-app}

In this subsection our goal is to derive \cite[Theorem 3.4]{MQ-sub} and then the geometric applications \cite[Theorem 4.2]{MQ-sub} 
based on Theorem \ref{Thm:theorem-7} and Theorem \ref{Thm:capacity-estimate} by the same ideas as in \cite{AH, MQ-2025, MQ-survey, CY}, which 
would not work without Theorem \ref{Thm:capacity-estimate}.

\begin{theorem} (\cite[Theorem 3.4]{MQ-sub}) \label{Thm:main-potential}
Let $\G=(\R^N, \circ, \delta_\lambda)$ be a polarizable Carnot group.
And let $R^{\mu, D}_{\Delta_\G}$ be the potential of a nonnegative Radon measure $\mu$ supported in $D\subset \R^N$ for the sub-Laplacian 
$\Delta_\G$ in $\G$.  Suppose that  $S$ is a compact subset in $D$ and that
\begin{equation}\label{Equ:complete-carnot-intr}
\lim_{t\to 0} R^{\mu, D}_{\Delta_\G}(\gamma(t)) = \infty \text{ and } \int_0^\delta (R^{\mu, D}_{\Delta_\G} (\gamma(t)))^\frac 2{Q-2} dt = \infty,
\end{equation}
for any curve $\gamma(t): [0, \delta]\to \R^N$ with $\gamma(0)\in S$ and $\gamma\in C[0, \delta]\bigcap C^1(0, \delta]$
that is horizontal with respect to the sub-Riemannian structure of $\G$ and
some $\delta >0$, where $Q$ is the homogeneous dimension of $\G$. Then 
\begin{equation}\label{Equ:intr-dim-p}
dim_{\mathcal{H}_{\Delta_\G}} (S) \leq \frac {Q-2}2.
\end{equation}
\end{theorem}

\begin{proof} Our proof here is different from that of \cite[Theorem 3.4]{MQ-sub} and relies on Theorem \ref{Thm:theorem-7} and Theorem \ref{Thm:capacity-estimate} in this 
paper.
By Theorem \ref{Thm:capacity-estimate}, one is able to pick up a horizontal radial curve that escapes from the thin set $E$ in Theorem \ref{Thm:theorem-7} and make use
of \eqref{Equ:estimate-upper} to reach to a contradiction to \eqref{Equ:complete-carnot-intr} if assuming otherwise
$
dim_{\mathcal{H}_{\Delta_{\mathsf{N}}}} (S)> \frac {Q-2}2.
$
\end{proof}

Consequently, due to the fact that the nilpotent part of the Iwasawa decomposition of rank one simple groups is H-type (cf. \cite{CDKR-91}) and therefore polarizable
(cf. \cite[Proposition 5.6]{BT}), we then have

\begin{theorem}\label{Thm:main-geometry} (\cite[Theorem 4.2]{MQ-sub}) Assume $\mathsf{N}$ 
is the nilpotent part of the Iwasawa decomposition of rank one groups: 
$SU(n+1, 1)$, $Sp(n+1, 1)$ and $\mathbb{F}_{4(-20)}$. 
Let $S$ be a compact subset of a bounded and open subset $\Omega$ in $\mathsf{N}$. 
And let $\theta = u^\frac 4{Q-2}\theta_0$ be a 
conformal contact form on $\Omega\setminus S$, where $\theta_0$ is the standard contact form on $\mathsf{N}$ and $Q$ is the homogeneous 
dimension of $\mathsf{N}$. 
Suppose that the sub-Riemannian metric associated with $\theta$ is 
geodesically complete near $S$. Then
\begin{equation}
dim_{\mathcal{H}_{\Delta_{\mathsf{N}}}} (S) \leq \frac {Q-2}2
\end{equation}
when the scalar curvature $R_\theta$ of the conformal contact form $\theta$ is nonnegative.
\end{theorem}

\begin{proof} The proof uses the same idea from that in \cite{MQ-2025, MQ-survey, CY}. 
The argument to make use of potentials and Theorem \ref{Thm:main-potential} for singular solutions to
scalar curvature equations is the same as that presented in \cite[Section 4.2]{MQ-sub}. \eqref{Equ:complete-carnot-intr} is implied by the geodesic completeness
of the sub-Riemannian structure associated with $\theta = u^\frac 4{Q-2} \theta_0$.
\end{proof}

With some efforts as in \cite{MQ-2025}, Theorem \ref{Thm:main-geometry} holds even when the scalar curvature is not pointwisely nonnegative analogous to that in
\cite[Theorem 1.3]{MQ-2025}. A version of Theorem \ref{Thm:main-geometry} for CR geometry has appeared in \cite[Theorem 1.1]{CY}. 

%We notice that, even
%the right-hand side of the equations with singular solutions is an integrable function away from the singularities, the right-hand side of the equations with singular
%solutions in distributional sense is a Radon measure which in general is not absolutely continuous.

%%%%%%%%%%%%%%%%%%%%%%%%%%%%%%%%%%%%%%%%%%
%%%%%%%%%%%%%%%%%%%%%%%%%%%%%%%%%%%%%%%%%%

\vskip0.1in
\noindent$\mbox{}^\dag$School of Mathematical Science and LPMC, Nankai University, Tianjin, China; \\e-mail: 
msgdyx8741@nankai.edu.cn 
\vspace{0.2cm}

\noindent $\mbox{}^\ddag$ Department of Mathematics, University of California, Santa Cruz, CA 95064; \\
e-mail: qing@ucsc.edu

\end{document}